\documentstyle[11pt]{article}
\begin{document}
\hoffset-2cm \oddsidemargin 1.0in
\title{{\Large \bf Loxodromic M\"{o}bius Transformations with Disjoint Axes}}
\author{\normalsize  Tan-Ran Zhang }
\date{}
\maketitle \baselineskip 18pt
\begin{minipage}{130mm}
\begin{abstract}
{This paper is concerned with Loxodromic M\"{o}bius Transformations
with disjoint Axes in a Kleinian Group. We study mainly the distance
between their axes, and give some estimates about their translation
lengths. }
\end{abstract}
\end{minipage}\\
\\
\\
\textbf{{\bf \large \textbf{1. Introduction}}}\\
\\
Let the set $$\textit{\textbf{H}} ^{3}=\{{(x_{1}, x_{2},
x_{3})\in{\textit{\textbf{R}} ^{3}}: x_{3}>0}\}$$ be the $hyperbolic
3-space$, and its metric is the complete Riemannian metric
$ds=|dx|/x_{3}$. Let $\mathcal{M}$ denote the group of all
M\"{o}bius Transformations of extended complex plane
$\textit{\textbf{\^{C}}}=\textit{\textbf{C}}\bigcup\{\infty\}$. A
Fuchsian group $G$ is a discrete subgroup of $\mathcal{M}$ with an
invariant disc $D \in \textit{\textbf{\^{C}}}$. The group of
orientation preserving isometries is denoted by
$Isom^{+}(\textit{\textbf{H}} ^{3})$, then a Kleinian group $G$ is a
discrete nonelementary subgroup of $Isom^{+}(\textit{\textbf{H}}
^{3})$. So it is easy to see, this kind of Kleinian group $G$ is not
virtually abelian, and it is also can be regarded as the extension
of a nonelementary Fuchsian group acting in $\textit{\textbf{H}} ^{3}$. This paper is
concerned with Loxodromic M\"{o}bius Transformations with disjoint
Axes in a Kleinian group, and we will give some estimates about
their translation lengths. For each M\"{o}bius Transformation

$$f=\frac{az+b}{cz+d}\in{\mathcal{M}}, \, \;  ad-bc=1,$$
the $2\times 2$ complex matrix $$A= \left(
\begin{array}{cc}
a& b  \\
c& d
\end{array}
\right)\in{SL(2,\textit{\textbf{C}})}$$ \\
induces $f$, and set ${\textrm{tr}(f)}^{2}={\textrm{tr}(A)}^{2}$,
where $\textrm{tr}(A)=a+d$, denotes the trace of the matrix $A$. And
for each $f$ and $g$ in $\mathcal {M}$, the commutator $[f,g]$ of
$f$ and $g$ is $fgf^{-1}g^{-1}$. We call the three complex numbers
$$ \gamma(f,g)=\textrm{tr}([f, g])-2, \quad \beta(f)=\textrm{tr}^{2}(f)-4,
\quad \beta(g)=\textrm{tr}^{2}(g)-4$$ the parameters of the
two-generator group $<f,g>$. These parameters are independent of the
choice of matrix representations for $f$ and $g$ in
$SL(2,\textit{\textbf{C}})$,   and they determine
uniquely up to conjugacy whenever $\gamma(f,g)\neq0$.\\

Let $f \in {\mathcal{M}}$ not be the identity, then \\
(1)\, $f$ is parabolic if and only if $\beta(f)=0$, and then $f$ is
conjugate to $z\rightarrow z+1$.\\
(2)\, $f$ is elliptic if and only if $\beta(f)\in{[-4,0)}$, and then
$f$ is conjugate to
$z\rightarrow\mu z$, where $|\mu|=1$.\\
(3)\, $f$ is loxodromic if and only if $\beta(f)\not\in{[-4,0]}$,
and then $f$ is conjugate to $z\rightarrow\mu z$, where $|\mu|>1$;
and $f$ is hyperbolic if $\mu>0$, $f$ is strictly loxodromic if
$\mu<0$ or $\mu$ is not real.\, Furthermore,\,
$\beta(f)=\mu-2+\mu^{-1}$.\\

\par A parabolic or hyperbolic element $g$ of a Fuchsian group $G$ is said to be $primitive$ if and only if
$g$ generates the stabilizer of each of its fixed points. If $g$ is elliptic, it is $primitive$ when it
generates the stabilizer and has an angle of rotation of the form $2\pi/n$.
And each  M\"{o}bius transformation of $\textit{\textbf{\^{C}}}=\partial \textit{\textbf{H}}
^{3}$ has a natural extension uniquely via the Poincar\'{e}'s way
 to an orientation-preserving isometry of hyperbolic 3-space
$\textit{\textbf{H}}^{3}$, see $[1]$. Then  Kleinian groups equal to
discrete M\"{o}bius groups. \\
\par If $f \in{\mathcal{M}}$ is nonparabolic, then $f$ has two fixed
points in $\textit{\textbf{\^{C}}}$ and the hyperbolic line
$(geodesic)$ joining these two fixed points in $\textit{\textbf{H}}$
is called the axis of $f$, denoted by $\mathrm{ax}$$(f)$. In this
case $f$ translates along ${\mathrm{ax}}(f)$ by an amount
${\mathrm{t}}$$(f)\geq0$, and ${\mathrm{t}}$$(f)$ is called the
translation length of $f$. $f$ rotates about ${\mathrm{ax}}(f)$ by
an angle $\theta(f)\in{(-\pi,\pi]}$, and
$$\beta(f)=4\sinh^{2}\left(\frac{{{\mathrm{t}}}(f)+i\theta(f)}{2} \right).\eqno  (1.1)$$
It then follows from (1.1) that\\
$$
\cosh({\mathrm{t}}(f))=\frac{|\beta(f)+4|+|\beta(f)|}{4}$$ and
$$
\cosh(\theta(f))=\frac{|\beta(f)+4|-|\beta(f)|}{4}$$
(cf.\,(15), (17) and (18) in [3]).\\
\par If $f, g \in{\mathcal{M}}$ are nonparabolic and if $\alpha$ is the
hyperbolic line in $\textit{\textbf{H}} ^{3}$ that is orthogonal to
${\mathrm{ax}}(f)$ and ${\mathrm{ax}}(g)$,   then
$$\frac{4\gamma(f,g)}{\beta(f)\beta(g)}=\sinh^{2}(\delta+i\phi),\eqno(1.2)$$
where $\delta$ is the hyperbolic distance between ${\mathrm{ax}}(f)$
and ${\mathrm{ax}}(g)$ and $\phi\in{[0, \pi]}$ is the angle between
the hyperplanes in $\textit{\textbf{H}} ^{3}$ that contain
${\mathrm{ax}}(f)\cup\alpha$ and ${\mathrm{ax}}(g)\cup\alpha$
respectively (see Lemma 4.2 in [4]). In particular if
${\mathrm{ax}}(f)$ and ${\mathrm{ax}}(g)$ are in the same hyperplane
then
$$\frac{4\gamma(f,g)}{\beta(f)\beta(g)}=\sinh^{2}(\delta).\eqno  (1.3)$$
\par Finally, there is a definition of triangle groups. A group $G$ of isometries of the hyperbolic plane
 is said to be of $type (\alpha, \beta, \gamma)$ if and only if $g$
 is generated by the reflections across the sides of some triangle
 with angles $\alpha, \beta$ and $\gamma$. A group $G$ is a  $(p, q,
 r)-Triangle$ $group$ if and only if $G$ is a conformal group of type
 $(\pi/p, \pi/q, \pi/r)$. We call $G$ a $Triangle$ $group$ if it is a $(p, q,
 r)-$ Triangle group for some integers $p, q$ and $r$. Triangle
 groups is an important class of Fuchsian groups. Roughly speaking,
 these are the discrete groups with the more closely packed orbits
 and the smallest fundamental regions, especially the (2, 3, 7) triangle
group, which is a simple but powerful example for kinds of extremal
conditions. It also can deduce the following two numbers that occur
frequently in this paper:
$$c=2(\cos(2\pi/7)+\cos(\pi/7)-1)=1.048...,$$
$$d=2(1-\cos(\pi/7))=0.198....$$\\
\par F. W. Gehring and G. J. Martin have given some similar results,
see [2] in details. Their research is mainly in the condition of
intersecting axes. There are two Theorems form theirs.\\
\textbf{{\large T\small HEOREM \large A}.} \quad \textit{ If $<f,g>$
is a Kleinian group, if $f$ and $g$ are hyperbolics, and if
${\mathrm{ax}}(f)$ and ${\mathrm{ax}}(g)$ intersect at an angle
$0<\phi<\pi$, then
$$
\sinh({\mathrm{t}}(f)/2)\sinh({\mathrm{t}}(g)/2)\sin(\phi)\geq
\lambda^{2},$$ where $\lambda=0.471...$\,. The constant $\lambda$ is
sharp and the exponent of $\sin(\phi)$ can not be replaced by a
constant greater than 1. }\\
\textbf{{\large T\small HEOREM \large B}.} \quad \textit{ If $<f,g>$
is a Kleinian group, if $f$ and $g$ are loxodromics with axes that
intersect at an angle $0<\phi<\pi$, then
$$
\mid \beta(f)\beta(g) \mid \sin^{4/3} (\phi) \geq b,$$ where
$0.777\geq b \geq 0.884$. }
\\
\par The main theorems of this paper are follows.\\
\textbf{{\large T\small HEOREM \large 1}.} \quad \textit{If $<f,g>$
is a Kleinian group,  $f$ and $g$ are hyperbolic, and if
\\${\mathrm{ax}}(f)$ and ${\mathrm{ax}}(g)$ do not intersect
then$$\sinh({\mathrm{t}}(f)/2)\sinh({\mathrm{t}}(g)/2)\sinh(\delta)\geq\sqrt{d}/2,\eqno
(1.4)$$ where $ \delta $ is the  distance between
${\mathrm{ax}}(f)$ and ${\mathrm{ax}}(g)$.}\\
\textbf{{\large T\small HEOREM \large 2}.} \quad \textit{If $<f,g>$
is discrete, $f$ and $g$ are loxodromics, if ${\mathrm{ax}}(f)$ and
${\mathrm{ax}}(g)$ do not intersect, and if
$\sinh(\delta)\leq1$,
then$$|\beta(f)\beta(g)|\sinh^{4/3}(\delta)\geq 4d ,\eqno(1.5)$$
where $\delta$ is the hyperbolic distance between ${\mathrm{ax}}(f)$
and ${\mathrm{ax}}(g)$.}\\
\textbf{{\large T\small HEOREM \large 3}.} \quad \textit{If $<f,g>$
is discrete, if $f$ and $g$ are loxodromics, ${\mathrm{ax}}(f)$ and
${\mathrm{ax}}(g)$ do not intersect, and
$\sinh(\delta)\leq1$, then $$
\sinh({\mathrm{t}}(f))\sinh({\mathrm{t}}(g))\sinh^{4/3}(\delta)\geq
\frac{3b}{16\pi^{2}}, \eqno(1.6)$$ where $\delta$ is the hyperbolic
distance between ${\mathrm{ax}}(f)$ and ${\mathrm{ax}}(g). $}\\
\textbf{{\large T\small HEOREM \large 4}.} \quad \textit{If $<f,g>$
is discrete, if $f$ and $g$ are loxodromics with
$\beta(f)=\beta(g)$, and if ${\mathrm{ax}}(f)$ and
${\mathrm{ax}}(g)$ do not intersect. then
$$ \sinh({\mathrm{t}}(f))\sinh(\delta)\geq\frac{\sqrt{3d}}{2\pi},\eqno(1.7)$$
where $\delta$ is the hyperbolic distance between ${\mathrm{ax}}(f)$
and ${\mathrm{ax}}(g). $}\\
\textit{ \textbf{{\large T\small HEOREM \large 5}.} \quad If $<f,g>$
is discrete, if $f$ is loxodromic and $g$ is elliptic of order
$n\geq 3$,  and if ${\mathrm{ax}}(f)$ and ${\mathrm{ax}}(g)$ do not intersect,
then$$\sinh({\mathrm{t}}(f))\sin^{2}(\pi/n)\sinh^{2}(\delta)\geq\sqrt{3}a(n)/4\pi,
\eqno(1.8)$$
where $\delta$ is the hyperbolic distance between ${\mathrm{ax}}(f)$
and ${\mathrm{ax}}(g)$. When $n\geq5$, we also have
$$\sinh({\mathrm{t}}(f))\sin^{2}(\pi/n)\sinh^{2}(\delta)\geq\sqrt{3}\cos(2\pi/n)/2\pi.$$ }
\\
\\
\\
\\
\\
 \textbf{{\bf \large \textbf{2. Proofs of Theorems}}}
 \\
 \\
\par These lemmas are to be used.\\
\textbf{{\large L\small EMMA \large 1}.} \quad  \textit{If $<f,g>$
is a Fuchsian group, then} \textit{$$ |\gamma(f,g)|\geq d$$} (see [5]).\\

\par The following result is established in [7] and [8], and then is
 sharpened by Cao in [9]. \\
\textbf{{\large L\small EMMA \large 2}.} \quad \textit{If $<f,g>$ is
a Kleinian group and if either
$$ |\beta(f)|\leq c \; \; or \; \; \beta(f)=\beta(g),\eqno(2.1)$$
then $$|\gamma(f,g)|\geq d.\eqno(2.2)$$ This result is sharp under
either assumption in $(2.1)$.}
\\
\textbf{{\large L\small EMMA \large 3}.} \quad \textit{For each
loxdromic M\"{o}bius transformation $f$ there exists an integer
$m\geq 1$ such
that$$|\beta(f^{m})|\leq\frac{4\pi}{\sqrt{3}}\sinh({\mathrm{t}}(f))
\eqno(2.3)$$} (see [10]).

\noindent\textbf{{\large L\small EMMA \large 4}.} \quad \textit{If
$<f,g>$ is a Kleinian group, $f$ is elliptic of order $n\geq 3$, and
$g$ is not of order 2, then $$|\gamma(f, g)|\geq a(n)$$ where
$$a(n)= \left\{\begin{array}{cc}
2\cos(2\pi/7)-1   & $if$  \;\ n=3,\\
\hspace*{-0.7cm} 2\cos(2\pi/5)     & \, \, \, \;  $if$  \; \ n=4,5,\\
\hspace*{-0.7cm} 2\cos(2\pi/6)     & $if$  \; \ n=6,\\
2\cos(2\pi/n)-1   & $if$  \; \ n\geq7
\end{array} \right.\eqno(2.4)
$$
\\}
(See [4]).\\

\textbf{{\large P\small ROOF OF THEOREM \large 1}.} Let
$S$ be the
hyperplane in $\textit{\textbf{H}} ^{3}$ determined by
${\mathrm{ax}}(f)$ and ${\mathrm{ax}}(g)$, then $S$ is invariant
under $G$, $F=G|S$ is conjugate to a Fuchsian group and
$$|\gamma(f,g)|\geq d \eqno(2.5)$$
by Lemma 1.
\par Next since $f$ and $g$ are hyperbolic, $\theta(f)=\theta(g)=0$
and
$$|\beta(f)|=4\sinh^{2}({\mathrm{t}}(f)/2),\,\; |\beta(g)|=4\sinh^{2}({\mathrm{t}}(g)/2)\eqno(2.6)$$
by (1.1). Thus
$$16\sinh^{2}({\mathrm{t}}(f)/2)\sinh^{2}({\mathrm{t}}(g)/2)\sinh^{2}(\delta))
=|\beta(f)||\beta(g)|\sinh^{2}(\delta)$$
$$\hspace{4.7cm}=4|\gamma(f, g)|$$
$$\hspace{3.8cm}\geq4d$$
by (1.3), (2.5) and (2.6), then we obtain (1.4).\quad $\Box$
\par An example due to J$\o$rgensen [6] shows that there exists no
absolute lower bound for $|\gamma(f,g)|$ when $<f,g>$ is a Kleinian
group.\\

\noindent \textbf{{\large P\small ROOF OF THEOREM \large 2}.} By
(1.2),
$$|\beta(f)\beta(g)|\sinh^{2}(\delta)=|4\gamma(f,g)|. \eqno(2.7)$$\\
We want to find a lower bound for
$$u=|\beta(f)\beta(g)|\sinh^{4/3}(\delta).$$
By relabeling, we may assume that $|\beta(f)|\leq|\beta(g)|.$
\par If $|\beta(f)|\leq c=1.048...$, then $\gamma(f,g)\neq
0$, by (2.7), and $<f,g>$ is a Kleinian group. Thus $|\gamma(f,
g)|\geq d$ by Lemma 2, and $$u\geq
|\beta(f)\beta(g)|\sinh^{2}(\delta)=|4\gamma(f,g)|\geq
4d=0.792...\,.$$ Next, if $|\beta(f)|\geq c$ , then
\begin{eqnarray}
      c^{-2}u^{3}+4u^{3/2} & = & |\beta(f)\beta(g)|^{3}c^{-2}\sinh^{4}(\delta)+4|\beta(f)\beta(g)|^{3/2}
      \sinh^{2}(\delta)\nonumber\\
                           & \geq & |\beta(f)\beta(g)|^{2}\sinh^{4}(\delta)+4|\beta(f)
                           \beta(g)|\sinh^{2}(\delta)|
                           \beta(f)|\nonumber\\
                           & \geq & 16|\gamma(f,g)|^{2}+16|\gamma(f,g)||\beta(f)|
                           \nonumber\\
                           & \geq & 16dc=3.320...\nonumber\end{eqnarray}
by (2.8) and Lemma 2, and we obtain$$u>0.798...\,.$$
Thus (1.5) follows. \quad $\Box$\\

\par The fact that $\sinh(\delta)\leq1$ is necessary. When
$\sinh(\delta)>1$, $u$ has no uniform lower bound. There is a
counterexample. For real $\lambda$, $\mu$ and $\delta$, define$$f=
\left(
\begin{array}{cc}
\ \cosh(\lambda)& e^{\delta}\sinh(\lambda)  \\
e^{-\delta}\sinh(\lambda)& \ \cosh(\lambda)
\end{array}
\right), \quad  g= \left(
\begin{array}{cc}
\ \cosh(\mu)& \sinh(\mu)  \\
\sinh(\mu)& \ \cosh(\mu)
\end{array}
\right). $$ It is clear that $f$ and $g$ are both loxodromic. The
axis of $f$ has endpoints $\pm e^{\delta}$ and the axis of $g$ has
endpoints $\pm 1$. By Reference [11], their axes are disjoint and
the distance between the two axes is $\delta$. Moreover,
$\beta(f)=4\sinh^{2}(\lambda)$ and $\beta(g)=4\sinh^{2}(\mu)$.
\par We suppose that $\textrm{tr}(fg^{-1})=-2$, so $fg^{-1}$ is parabolic
in $PSL(2,\textit{\textbf{C}}).$ For all $\lambda>0$ and all $\mu>0$
the group $<f,g>$ is clearly discrete and non-elementary. Indeed it
is a free group. Moreover,
$$-2=\textrm{tr}(fg^{-1})=2-2\sinh(\delta)\sinh(\lambda)\sinh(\mu).$$
Therefore, $$\sinh(\delta)=2/(\sinh(\lambda)\sinh(\mu)).$$ We have
\begin{eqnarray}
      u & = & |\beta(f)\beta(g)|\sinh^{4/3}(\delta)\nonumber\\
                           & = & 16\sinh^{2}(\lambda)\sinh^{2}(\mu)(2/(\sinh(\lambda)
                           \sinh(\mu)))^{4/3}\nonumber\\
                           & = &
                           2^{16/3}\sinh^{2/3}(\lambda)\sinh^{2/3}(\mu).
                           \nonumber\end{eqnarray}
Fixing $\mu$ and letting $\lambda$ tend to zero shows that there is
no uniform lower bound on $u$.\\

\noindent \textbf{{\large P\small ROOF OF THEOREM \large 3}.} By
Lemma 3 we can choose integers $m, \ n\geq 1 $ such
that$$|\beta(f^{m})|\leq
\frac{4\pi}{\sqrt{3}}\sinh({\mathrm{t}}(f)),\quad |\beta(g^{n})|\leq
\frac{4\pi}{\sqrt{3}}\sinh({\mathrm{t}}(g)).\eqno(2.8)$$ Then
$<f^{m},g^{n}>$ is Kleinian and we obtain
$$\left(\frac{4\pi}{\sqrt{3}}\right)^{2}\sinh({\mathrm{t}}(f))\sinh({\mathrm{t}}(g))
\sinh^{4/3}(\delta)\geq |\beta(f^{m})\beta(g^{n})|\sinh^{4/3}(\delta)\geq
b$$ by (2.8) and (1.5). This implies (1.6). \quad $\Box$
\par We see immediately there is not a lower bound for
$\max({\mathrm{t}}(f),{\mathrm{t}}(g))$.\\

\noindent \textbf{{\large P\small ROOF OF THEOREM \large 4}.} By
Lemma 3 we can choose an integer $m\geq1$ such that
$$|\beta(f^{m})|\leq \frac{4\pi}{\sqrt{3}}\sinh({\mathrm{t}}(f)).$$
Then $<f^{m},g^{n}>$ is Kleinian with $\beta(f^{m})=\beta(g^{m})$
and we obtain
\begin{eqnarray}
\frac{4\pi}{\sqrt{3}}\sinh({\mathrm{t}}(f))\sinh(\delta)& \geq & (|\beta(f^{m})\beta(g^{m})|
\sinh^{2}(\delta))^{1/2}\nonumber\\
                                           & = &
                                           (4|\gamma(f^{m},g^{m})|)^{1/2}\geq
                                           2\sqrt{d}\nonumber\end{eqnarray}
from (2.3), (1.3), and (2.2). This implies (1.7).\quad $\Box$\\

\noindent \textbf{{\large P\small ROOF OF THEOREM \large 5}.} We may
assume without loss of generality that $g$ is a primitive elliptic.
Next, by Lemma 3 we can choose an integer $m\geq1$ such that
$$\beta(f^{m})\leq \frac{4\pi}{\sqrt{3}}\sinh({\mathrm{t}}(f)).$$
Then $<f^{m},g>$ is Kleinian and $f^{m}$ is not of order 2, so we
obtain$$\frac{4\pi}{\sqrt{3}}\sinh({\mathrm{t}}(f))\sin^{2}(\pi/n)\sinh^{2}(\delta)\geq
|\gamma(f^{m},g)|\geq a(n)\eqno(2.9)$$ from (1.3) and Lemma 4, where
$a(n)$ is as in (2.4). Then (2.9) yields (1.8). \quad $\Box$
\\
\\
\\
 \textbf{{\bf \large \textbf{3.Remarks}}}
 \\
 \\
\par For elliptic transformations with disjoint axes whenever intersect,
there are also similar results as Theorem 1 and Theorem 2. Theorem A
holds with equality if $f$ and $g$ are hyperbolic generators for the
$(2, 3, 7)-$ triangle group with
$$ par(<f, g>)=(-d, c, c)=(-0.198...,1.048...,1.048...).$$ \
\\
\\
\\
\\
\\
\\
\\
\\
\\

\end{document}